\documentclass[12pt]{article}
\usepackage{amsmath,amsthm,amstext,amsfonts,amssymb,amscd}
\usepackage{epsfig,psfrag}
\usepackage{epstopdf}
\usepackage{makeidx}
\usepackage[latin1]{inputenc}

\catcode`\@=11

\def\ps@headings{\let\@mkboth\markboth
    \def\@oddfoot{}\def\@evenfoot{}
    \def\@evenhead{{\rm \thepage}\hfil { \rm \leftmark}}
    \def\@oddhead{{\rm \rightmark}\hfil {\rm \thepage}}
    \def\chaptermark##1{\markboth{\ifnum \c@secnumdepth >\m@ne
          \@chapapp\ \thechapter. \ \fi ##1}{}}%
    \def\sectionmark##1{\markright{\ifnum \c@secnumdepth >\z@
       \thesection. \ \fi ##1}}}

\ps@headings

\catcode`\@=12

\newtheorem{proposition}{Proposition}[section]
\newtheorem{definition}[proposition]{Definition}
\newtheorem{lemma}[proposition]{Lemma}
\newtheorem{theorem}[proposition]{Theorem}
\newtheorem{cor}[proposition]{Corollary}

\newcommand {\CP}{{\cal P}}

\newcommand {\CT}{{\cal T}}

\newcommand{\N}{{\mathbb N}}

\newcommand{\R}{{\mathbb R}}

\newcommand{\Z}{{\mathbb Z}}

\newcommand{\eps}{\varepsilon}
\newcommand{\rro}{\varrho}

\newcommand{\8}{\infty}
\newcommand{\disp}{\displaystyle}
\newcommand{\ninf}{{n\rightarrow +\8}}

\renewcommand{\S}{{\Sigma}}

\theoremstyle{definition}
\newtheorem{remark}{Remark}

\begin{document}

\title{Equilibrium States for Partially Hyperbolic Horseshoes}
\author{R. Leplaideur, K. Oliveira, and I. Rios}

\maketitle

\begin{abstract}
In this paper, we study ergodic features of invariant measures for
the partially hyperbolic horseshoe at the boundary of uniformly hyperbolic diffeomorphisms
constructed in \cite{DHRS07}. Despite the fact that the non-wandering set is a horseshoe, it contains intervals. We prove that every recurrent point has non-zero Lyapunov exponents and all
ergodic invariant measures are hyperbolic. As a consequence, we obtain  the existence of equilibrium
measures for any continuous potential. We also obtain an example
of a family of $C^\infty$ potentials with phase transition.
\end{abstract}


\section{Introduction}

Smale's horseshoes and  geodesic flows in negative curved surfaces are, doubtless,   landmarks in the study of dynamical systems theory. These examples are source for the general theory of \emph{hyperbolic maps},  extensively studied from the structural and ergodic points of view. From the structural point of view, they are stable and open with respect to appropriate  $C^1$ topology, exhibiting absolutely continuous invariant dynamical manifolds and they are sensitive with respect to small perturbations in initial conditions. From the ergodic point of view, there exist equilibrium measures for any continuous potential and every  H\"older potentials admits an unique equilibrium state (on each transitive component). These measures are Gibbs
states of the system, have full support, and their microscopic features are now very well  understood.

To extend this  theory beyond the hyperbolic setting is, nowadays, a very challenging task. In this direction, a very successful concept that extends the notion hyperbolicity, allowing some of its main consequences to be achieved, is the notion of \emph{partial hyperbolicity}. This notion is a weak version of hyperbolicity that preserves many of its features. Many authors have successfully  established results for partial hyperbolic systems: existence of absolutely continuous invariant manifolds (\cite{BP74}),  robust transitivity and generic properties (\cite{DPU99,BDP03, BDPR03}), ergodic properties, as existence of SRB measures and stable ergodicity (\cite{BoV00,ABV00,BDP02}).

Concerning equilibrium states, some  fruitful approaches  by many authors have established existence or uniqueness beyond the hyperbolic setting, when the system under consideration has a specific structure. For instance, for interval maps, rational functions of the sphere, and H\'enon-like
maps we cite \cite{BrK98,DU92,WY01}; for countable Markov shifts and piecewise expanding maps, \cite{BS03, Sa03,Yu03}; for horseshoes with tangencies at the boundary of hyperbolic
systems,  \cite{LR1}; for  higher dimensional local diffeomorphisms, \cite{Ol03,AMO04, OV06}, just to mention a few of the most recent works.
 Philosophically, a few restrictions can be imposed to the system to provide that all candidates for equilibrium measures have their  exponents bounded  away from zero. This gives to these non-uniformly hyperbolic maps the ``flavor" of uniform hyperbolicity.

 In this work we deal with a family of three dimensional partially hyperbolic horseshoes  $F$ at the boundary of uniformly hyperbolic diffeomorphisms. Each one of these horseshoes displays  a heteroclinic cycles and its non-wandering set $\Lambda$ contains intervals. They were constructed in \cite{DHRS07} as time zero of  bifurcations of families $F_t$ of partially hyperbolic horseshoes. Here, we give a complete description of Lyapunov exponents in the  central direction for ergodic measures, and prove that they are hyperbolic. As consequence of this, we get that \emph{any} continuous potential admits equilibrium states.

   Concerning uniqueness, we prove that the family $\phi_t=t\log|DF|_{E^c}|$ has a phase transition: there exists a $t_0>0$ such that $\phi_{t_0}$ admits at least two different equilibrium states.  In view of the recent results of \cite{OV07}, it is likely that there is a unique equilibrium state for any $t$ small enough. In fact, this seems to be true for any H\"older potential $\phi$ such that $\sup \phi - \inf \phi$ is smaller than some constant that depends only on the topological entropy and the expansion/contraction  rates of $F$.

\textbf{Acknowledgments:} We are grateful to F. Abdenur, L. D\'iaz and M. Viana for useful conversations.
Most of this work was carried out at Universit\'e de Bretagne Occidentalle (UBO) at Brest, France. K.O. is also thankful to  Penn State University for the hospitality during the final preparation of this manuscript.
This paper was partially supported by
CNPq, CAPES, Faperj, and UBO.

\subsection{Definition of the family of diffeomorphisms}\label{s.model}

 We consider in $\R^3$ a
family of horseshoe maps
$F=F(\lambda_0,\lambda_1,\beta_0,\sigma,\beta_1) \colon R\rightarrow
\R^3$, on the cube $R=I^3$, where $I$ denote the interval $I=[0,1]$.
Define the sub-cubes $R_0=I\times I\times [0,1/6]$, and
$R_1=I\times I\times [5/6,1]$ of $R$. The restrictions $F_{i}$ of $F$
to $R_i$, $i=0,1$, are defined by:

\begin{itemize}
\item
$F_{0}(x,y,z)=F_0(x,y,z)=(\lambda_0 x, f(y), \beta_0 z)$,
with $0<\lambda_0<1/3$, $\beta_0>6$ and $f$ is the time one map of a
vector field to be defined later;
\item
$F_{1}(x,y,z)=(3/4-\lambda_1 x, \sigma(1-y), \beta_1 (z-5/6))$,
with $0<\lambda_1<1/3$, $0<\sigma<1/3$ and $3<\beta_1<4$.
\end{itemize}

\begin{figure}[htb]

\begin{center}

\psfrag{f}{$f$}
\psfrag{0}{$0$}
\psfrag{1}{$1$}
\includegraphics[height=1.2in]{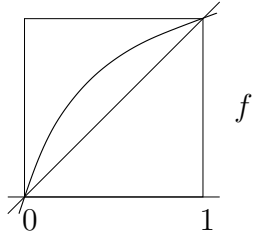}
 \caption{\label{fig.mapf}The central map $f$}

\end{center}

\end{figure}

The map $f:I\to I$ is defined as the time one of  the vector field
$$
x'=-x(1-x).
$$
This map is depicted in Figure~\ref{fig.mapf}.

Observe that $f(0)=0$ and $f(1)=1$, and , for every $y\neq 0$,
\begin{equation}
\label{e.fn}
f^n(y)=\frac{1}{1-\left( 1-\frac{1}{y}\right) e^{-n}}.
\end{equation}
We also have
\begin{equation}
\label{e.dfn}
\displaystyle{(f^n)'(y)=\frac{-\frac{e^{-n}}{y^2}}{\left( 1-\left(
1-\frac{1}{y}\right) e^{-n}\right) ^2}= -\frac{e^{-n}}{y^2}\left(
f^n(y)\right) ^2.}
\end{equation}
Note that $f'(0)= e$ and $f'(1)=1/e$. Since we have  $f(0)=0$ and $f(1)=1$,
 the point $Q=(0,0,0)$ is a fixed
saddle of index $1$ of $F$, and the point
 $P=(0,1,0)$ is a fixed saddle of index $2$ of $F$.

In \cite{DHRS07}, the authors proved that the diffeomorphism $F$ is simultaneously at the boundary of the sets of uniformly hyperbolic systems and  robustly non-hyperbolic systems.
In fact, it is reached as the first bifurcation of a one-parameter family of $C^{\infty}$ diffeomorphisms of the space $\R^3$, whose unfolding leads to robust non-hyperbolic behavior.  Here we state some other properties of the diffeomorphism $F$, see \cite{DHRS07} for proofs.

\begin{enumerate}
\item
The diffeomorphism $F$ has a heterodimensional cycle associated to the
saddles $P$ and $Q$.
\item
The homoclinic class of $Q$ is trivial and
the homoclinic class of $P$ is non-trivial and contains the saddle
$Q$, thus $H(Q,F)$ is properly contained in $H(P,F)$.
\item
There is a
surjection
$$
\Pi\colon H(P,F)\to \Sigma_{11}, \qquad \mbox{with} \quad \Pi \circ F=
\sigma \circ \Pi,
$$
and infinitely many central curves $C$ such that every $C$ contains
infinitely many points of the homoclinic class of $H(P,F)$ and
$\Pi(x)=\Pi(y)$ for every pair of points $x,y\in C \cap H(P, F)$.
These intervals consist of non-wandering points.
\end{enumerate}

For the rest of the paper, $\Lambda$ will denote the maximal
invariant set in the cube. Namely
$$
\disp \Lambda=\bigcap_{n\in\Z}F^{-n}(R). $$ For $X=(x^s,x^c,x^u)$ in
$\Lambda$, we denote by $W^u(X)$ and $W^s(X)$ the strong unstable and
strong stable leaves of $X$. The central leaf $W^c(X)$, will denote
the set of points on the form $(x^s,y,x^u)$,
 with $y\in I$. Given an  ergodic invariant measure $\mu$ we define the \emph{central
Lyapunov exponent} as:
$$
\lambda_\mu^c=\int \log |DF|_{E^c}| \, d\mu.
$$ Note that, since $E^c$ is one dimensional and $\mu$ is ergodic,
$$
\lambda_\mu^c = \lim_\ninf \frac{1}{n}\log |DF^n(p)|_{E^{c}}| ,
$$ for $\mu$ almost every point $p\in \Lambda$.

\section{Statement of the main results}

Our first result is a description of central Lyapunov exponents for
any ergodic invariant measure. We prove that the central Lyapunov
exponents of these measures are negative, except for the measure
$\delta_Q$. Moreover, using this information, we are able to prove
the existence of equilibrium measures associated to any continuous
potential. 

\begin{theorem}\label{t.A} The following properties of $F$ hold
true:

\begin{enumerate}

\item For any recurrent point $p$ different from $Q$:
$$
\liminf_{\ninf} \frac{1}{n} \log |DF^n(p)|_{E^{c}}|\leq 0.
$$ Moreover, any ergodic invariant measure for $F$
different from $\delta_Q$ has negative central Lyapunov exponent.

\item If $\mu_n$ is a sequence of ergodic invariant measures such that $\lambda_\mu^c$
converges to zero, then $\mu_n$  converges to
$\frac{\delta_Q+\delta_P}{2}$ in the weak$^\star$ topology.

\end{enumerate}

\end{theorem}

Let $\phi: R\to\R$ be a continuous function. Let $\eta$ be a $F$-invariant probability measure. The $\eta$-pressure of the potential $\phi$ (or equivalently the $\phi$-pressure of the measure $\eta$) is defined by 
$$ h_\eta(f)+\int\phi\,d\eta.$$
We recall that $\eta$ is called  an \emph{equilibrium state}  for
the potential $\phi$ if its $\phi$-pressure maximizes the $\phi$-pressures
among all $F$-invariant probabilities. Our second result is:

\begin{theorem}\label{t.A'} Any continuous function  $\phi: R \rightarrow
\R$ admits an equilibrium state. Moreover, there exists a
residual set of $C^0(\Lambda)$ such that the equilibrium measure is
unique.
\end{theorem}

If $\mu$ is an equilibrium state for some continuous potential $\phi$, the $\phi$-pressure of $\mu$ is also the \emph{topological pressure} of $\phi$.
A natural  question  that arises from the previous theorem is if
H\"older regularity of $\phi$ implies uniqueness of the equilibrium
measure.   A negative
answer to this question for a particular potential is given in
Theorem~\ref{t.C} below.  Concerning this, we  are able to prove
that some restriction is necessary. We
prove that $\phi_t=t\log|DF|_{E^c}|$ admits a phase transition:

\begin{theorem}\label{t.C} Consider the one parameter family $\phi_t$ of $C^\infty$
potentials defined for $X=(x,y,z)\in R$ by
\begin{equation*}
\phi_t(X)=t\log|DF(X)|_{E^{c}}|=\left\{
                                        \begin{array}{ll}
                                          t\log f'(y), & \hbox{for $z\leq 1/6$;} \\
                                          t\log \sigma, & \hbox{for $z\geq 5/6$.}
                                        \end{array}
                                      \right.
\end{equation*}
There exists a positive real number $t_{0}$ such that:
\begin{enumerate}

\item  For $t>t_0$, $\delta_Q$ is the unique equilibrium state.

\item For $t< t_0$, any equilibrium state for
$\phi_t$ has negative central Lyapunov exponent. In particular, this measure
is singular with $\delta_Q$.

\item For $t= t_0$, $\delta_Q$ is an equilibrium state for
$\phi_t$, and there exists at least one other equilibrium state, singular with $\delta_Q$.

\end{enumerate}

\end{theorem}

\begin{remark}
In fact, $t_0$ can be defined as the supremum, among all $F$-invariant measures different from $\delta_Q$, of the expression $\left\{\frac{h_\mu(F)}{1-\lambda^c_\mu}\right\}$. 
Note that by Theorem \ref{t.A}, this number is well defined.
\end{remark}
\section{Central Lyapunov exponents}

In this section we study some interesting features of $F$. We
are able to prove that if $x \in \Lambda$ is recurrent and different from $P$ and $Q$, then $W^c(X)
\cap \Lambda = \{X\}$, despite the fact that $\Lambda$ contains central intervals. We also
prove that central Lyapunov exponents of any ergodic measure
different from $\delta_Q$ is negative. 

\subsection{Central Lyapunov exponents for recurrent points}\label{s.sif}

The main tool to prove the results in this section is the reduction
of the dynamics to a one-dimensional system of iterated functions.
Here we study these system, as well as some definitions and results
in \cite{DHRS07} that we need in this work.

Consider the maps $f_{0},f_{1}\colon I\to \R$ defined by
\begin{equation*}
\begin{array}{ll}
f_{0}(y)&=f(y),
\\
f_{1}(y)&=\sigma\, (1-y).
\end{array}
\end{equation*}

Given any $X=(x_0^s,x_0^c,x_0^u)\in \Lambda$ and $k\ge 0$, let $X_k=F^k(X)=(x_k^s,x^c_k,x^u_k)$. By the definition
of $F$, the central coordinate $x_k^c$,
of $X_k$ is
$$
x_{k}^c=
f_{i_{k-1}}\circ f_{i_{k-2}} \circ \cdots
\circ f_{i_0} (x_0^c),
$$
where the numbers $i_0,\dots,i_{k-1}\in\{ 0,1\}$ are determined by
the coordinates $x^u_{0},\dots,x^u_{k-1}$. In fact, the map $F$
admits a well defined projection along the central direction to
$I^2$, and this projection is conjugated to the shift in $\Sigma_{11}$.
Using this conjugacy, $F$ can be thought as the skew product

\begin{eqnarray*}
\tilde{F}\colon \Sigma_{11}\times I & \to & \Sigma_{11}\times I\\
                (\theta , x)  & \mapsto & (\sigma \theta, f_{\theta}(x)),
\end{eqnarray*}
where $f_{\theta}=f_{\theta_0}\in\{ f_0, f_1\}$.

In what follows, we consider the dynamics associated to the system of iterated functions (s.i.f.)
generated by $f_{0}$ and $f_{1}$, that we denote by $\mathfrak{F}$.

Given a sequence $(i_n)\in \Sigma_{11}^+$, for each given $k\ge 0$ we
consider the $k$-block $\rro_k=\rro_k(i_n)=(i_0,i_1,\dots,i_k)$
associated to $(i_n)$. For each $k$-block $\rro_k$, we consider the
map $\Phi_{\rro_k}$ defined by

$$
\Phi_{\rro_k}(x)= f_{i_k}\circ f_{i_{k-1}}\circ \cdots \circ f_{i_0}(x).
$$






The computation of the contraction in the central direction is based on an explicit computation of the derivative of the functions  $\Phi_{\rro_k}$. First, consider a point $y\in(0,1]$.
Then we have
$$
|(f_{1}\circ f^{\alpha}_{0})'(y)|= \left(\frac{w}{y\,
(1-y)}\right) \, \left( 1- \frac{w}{\sigma}\right), \quad
\mbox{where $f_0^\alpha(y)=1-w/\sigma$.}
$$
Note that $f_{1}\circ f^{\alpha}_{0}(y)=w$. This implies that, if we
chose a sequence  $(i_n')\in \Sigma_{11}^+$ such that $(i_n')$ is the
concatenation of blocks of type $(0, \dots ,0,1)$, with the 1's
occurring in the positions $k_i$, we have

\begin{equation}\label{e.equation}
\Phi_{\rro_{k_i}}^\prime (y) = \prod_{j=1}^{i}
 \frac{w_j \, (1-{w_j}/{\sigma})}{w_{j-1}\, (1-w_{j-1})} \quad
 \mbox{where $w_0=y$ and  $w_j=\Phi_{\rro_{k_j}}(y)$.}
\end{equation}

Observe that, if $w_j>0$, the factor of the product in 
(\ref{e.equation}) corresponding to it is strictly smaller than $1$.
Moreover, it is a decreasing function of $w_j\in[0,\sigma]$.

\begin{lemma}\label{s.contraction}
Let  $(i_n)\in \Sigma_{11}^+$ be a sequence with infinitely many $1$'s. Assume that $i_{0}=1$. Let $(n_{j})_{j\geq 0}$ be the sequence of positions of the $j+1$'s symbols $1$ for $(i_{n})$. Then, there exist a sequence of positive real numbers $(\delta_{j})_{j\geq 0}$ and a positive constant $C$ such that 
\begin{itemize}
\item[(i)] for every $y$ in $[0,1]$, 
$\disp
|\Phi_{\rro_{n_i}}^\prime (y)| 
 \leq C \prod_{j=1}^{i-1}
 \frac{
   1-{\delta_j}/{\sigma}}{ 1-\delta_{j}}$,
\item[(ii)] $C$ depends only on $n_{0}$,
\item[(iii)] each $\delta_{j}$ depends only on the $n_{i}$'s, $i\leq j$.
\end{itemize}
\end{lemma}

\begin{proof}
Let $\rro'$ be the block of $(i_n)$ starting at the first
symbol and finishing at the second $1$. Let $N =n_{0}$ be its size, and
$(i_n')$ be the sequence obtained from $(i_n)$ by removing $\rro'$.
Then, for $k>N$ and $y\in (0,1]$, 
\begin{equation}\label{equ1-rustiren}
\Phi_{\rro_{k}}^\prime
(y)=\Phi_{\rro_{(k-N)}'}^\prime
(\Phi_{\rro_{N}}(y)).\Phi_{\rro_{N}}^\prime (y).
\end{equation}
Let $A= \max \{
|\Phi_{\rro_{N}}^\prime (\xi)|,\xi\in I\}$. Note that $A$ only depends on $n_{0}$.

Let  $w_0=\Phi_{\rro_{N}}(y)$,  and
$w_j=\Phi_{\rro_{n_j-N}'}(w_0)$. Observe that $\Phi_{\rro_{N}}(I)\subset
(0,\sigma]$; we set 
$$\delta_0=\min\Phi_{\rro_{N}}(I) \mbox{ and }
\delta_j
=\min\Phi_{\rro_{n_j-N}'}([\delta_0,\sigma]).$$
Then, \eqref{e.equation} yields 
\begin{equation}\label{e.equation2-prime}
|\Phi_{\rro_{n_i-N}'}^\prime (w_0)| = \frac{w_i \, (1-{w_i}/{\sigma})}{w_0\, (1-w_0)} \prod_{j=1}^{i-1}
 \frac{ 1-{w_j}/{\sigma}}{1-w_{j}}
 \leq \frac{1}{3\delta_0(1-\delta_{0})} \prod_{j=1}^{i-1}
 \frac{
   1-{\delta_j}/{\sigma}}{ 1-\delta_{j}}.
\end{equation}
Therefore, \eqref {equ1-rustiren} and \eqref{e.equation2-prime} yield $(i)$, with $\disp C=\frac{A}{3\delta_0(1-\delta_{0})}$. Note that $C$ only depends on $n_{0}$. Moreover each $\delta_{j}$ only depends on the $n_{i}$'s, with $i\leq j$. This finishes the proof of the lemma.
\end{proof}

\begin{lemma}\label{lem-scholium2}
Let  $(i_n)\in \Sigma_{11}^+$ be a recurrent sequence for the shift such that $i_{0}=1$. Then there exist a real number $a$ in $(0,1)$ and an increasing sequence of times $(m_{j})_{j\geq 0}$ such that for every $y$ in $[0,1]$,
$$|\Phi_{\rro_{m_{j}}}'(y)|\leq C.a^{j},$$
where $C$ is obtained from $(i_{n})$ as in Lemma \ref{s.contraction}. 
\end{lemma}
\begin{proof}
Note that as the sequence $(i_{n})$ is recurrent, it has infinitely many symbols $1$. We can thus apply Lemma \ref{s.contraction}. In particular, we use the notations of its proof.

Since each factor in the product  in $(iii)$-Lemma \ref{s.contraction} is strictly less than  $1$, it remains to show that there are infinitely
 many factors bounded from above by a number strictly smaller than $1$. This is equivalent to show that there are infinitely
  many values of $j$ such that $\delta_j$ is uniformly bounded away from zero.

The first block of $\rro'$ is composed by $n_1-1$ zeros and one $1$.
This implies that $\Phi_{\rro_{n_1-N}'}[0,\sigma]\subset [f_1\circ
f_0^{n_1-1}(\sigma),\sigma]$, and so $\delta_1>f_1\circ
f_0^{n_1-1}(\sigma)$. By the recurrence of the sequence $(i'_{n})$,
this first block repeats infinitely many times. For each time $j$
that it repeats, using the same argument, we conclude that
$\delta_{j+1}\in [f_1\circ f_0^{n_1-1}(\sigma),\sigma]$. This
concludes the proof. 
\end{proof}

\begin{remark}
A direct consequence of Lemma \ref{lem-scholium2} is that any periodic point is hyperbolic, and if it is different from $Q$, it admits a negative Lyapunov exponent in the central direction.
\end{remark}

\begin{remark}\label{rem-scholmium2-smooth}
The hypothesis ``$(i_{n})$ recurrent'' is not necessary, and it can be replaced by the weaker assumption: ``One block of the form $(1,\underbrace{0\ldots,0}_{k},1)$, with a fixed $k$, appears infinitely many times in $(i_{n})$''. 
\end{remark}



\subsection{Proof of Theorem \ref{t.A}}

Let $X$ be a recurrent point for $F$ (for forward and backward iterations). Assume $X$ is different from $Q$ and $P$. Then $X$ is forward-recurrent for $F$. Let us consider the one-sided sequence $\Pi(X)^{+}$, which is recurrent and admits infinitely many symbols $1$. Hence,  we can apply Lemma \ref{lem-scholium2} to obtain 
$$
\liminf_{\ninf} \frac{1}{n} \log |DF^n(X)|_{E^{c}}|\leq 0.
$$ 
 This gives estimates for the forward iteration, but we can also get estimates for the backward iterations:

\begin{proposition}\label{prop-renaud}
Let $X$  be a recurrent point for $F$ (for forward and backward iterations) different from $Q$ and $P$. 
 Then
$$\cap_{n\in\Z}F^n(R)\cap W^{c}(X)=\{X\}.$$ 
\end{proposition}

\begin{proof}
Let  $\Pi(X)=(i_{n})\in\S_{11}$. 
Without loss of generality, we can assume that $i_0=1$. Let
$\rro_{k}$ denote any block $(i_0,\dots,i_k)$ of the sequence
$(i_n)$.  We denote by $(i_{n^+})$ the associated one-sided sequence. Again, we use vocabulary and notations from the proofs of Lemmas \ref{s.contraction} and \ref{lem-scholium2}. 

The infinite block $[(i_{0},i_{1},\ldots)]$ begins with the concatenation of the blocks $\rro_{n_{0}}$ and $\rro'_{n_{1}-n_{0}}$. The constant $C$ in Lemma \ref{s.contraction} only depends on $n_{0}$. The sequence  of $m_{j}$'s in Lemma \ref{lem-scholium2} is the sequence of appearances of the block $\rro'_{n_{1}-n_{0}}$ (``shifted'' to the end of the appearance).

By recurrence of $(i_n)$, we know that the block $\rro_{n_{0}}\rro'_{n_{1}-n_{0}}$ appears
infinitely many times in the sequence $(\dots, i_{-2}, i_{-1},
i_0)$. We consider   a decreasing sequence of integers $k_j\to -\infty$
such that $\sigma^{-k_j}((i_n))$ coincides with $(i_n)$ at the
positions $0,1,\dots,n_{1}$. We also ask that $k_j-k_{j+1}>n_{1}$.
 Lemma \ref{lem-scholium2} implies that for every $j$ and for every $y$ in $[0,1]$, 
\begin{equation}\label{equ1-correclem36}
|\Phi'_{[(i_{k_{j}},\ldots,i_{-1})]}(y)|\leq C. a^{j}.
\end{equation}

Let $L_j\subset I$ be the image of the interval $I$ by the map
$\Phi_{[(i_{n_j}, \dots ,i_{-1})]}$. Points
in $\cap_{n\in\Z}F^n(R)\cap W^{c}(X)$ have
their central coordinates  belonging to the intersection of the
sets $L_j$, $j>0$. Now,
\eqref{equ1-correclem36}
 implies that the
diameter of $L_j$ converges to zero. We also have that each $L_j$ is
non-empty, compact and $L_{j+1}\subset L_j$.
 Thus, their intersection is a single point. This completes the proof of the
 proposition.
\end{proof}

We define the \emph{cylinder} associated with the block
$\rro=(i_0,\dots,i_k)$ as follows:
$$
[\rro] = [i_0,\dots,i_k] = \{x \in \Lambda; F^j(x)\in R_{i_j},
\text{ for } j=0,\dots,k\}= \bigcap\limits_{j=0}^k F^{-j}(R_{i_j}) \cap \Lambda.
$$

The last  expression in the definition above  tells us that these sets
are always closed sets, since they are finite intersection of closed
sets $F^{-j}(R_{i_j}).$ We say that a point $p$ has \emph{positive
frequency} for a set $A \subset \Lambda$ if
$$
\gamma(p,A)=\liminf \frac{\# \{0\leq j<n;f^j(p)\in A\}}{n}>0.
$$
\begin{definition}
We say that a point $p$ is  of contractive type if 
$$\liminf_{\ninf}\frac1n\log|DF^n(p)_|{E^{c}}|<0.$$
\end{definition}

Next proposition is a tool to finish the proof of item \emph{1} in Theorem \ref{t.A}.

\begin{proposition}\label{prop.freq}
Let $p\in \Lambda$   be a point  with  positive frequency $\gamma>0$ at some
cylinder associated with a $l$-block $\theta=(0,0,\dots,0,1)$. Then it is of contractive type and its
central Lyapunov exponent is   less than a constant
$c(\gamma,l)<0$ that depends only on $\gamma$ and $l$.
\end{proposition}

\begin{proof}
We simply use Lemma \ref{lem-scholium2} and Remark \ref{rem-scholmium2-smooth}. There exist a constant $C=C(p)$ and $a\in(0,1)$, such that for every $n$ satisfying $F^{n}(p)\in \theta$, 
\begin{equation}\label{equ1-correcprop35}
|DF^{n+l}(p)|_{E^{c}}|\leq C(p)a^{\#\{0\leq j\leq n,\ F^{j}(p)\in A\}}.
\end{equation}
Note that $a$ depends only on the length of the cylinder $\theta$, hence on $l$. Moreover, \eqref{equ1-correcprop35} yields
$$\liminf_{\ninf}\frac1n\log|DF^{n+l}(p)|_{E^{c}}|\leq \gamma(p,\theta)\log a<0.$$
This finishes the proof of Proposition~\ref{prop.freq}.
\end{proof}

\begin{cor}\label{cor-lyapu-meas}
Every ergodic and $F$-invariant probability $\mu$ which is not
$\delta_Q$ has a negative Lyapunov exponent in the central
direction.
\end{cor}
\begin{proof}
First, since every point in $[0,0,\dots]\setminus\{Q\}$ is
attracted to $P$, the cylinder $[0,0,\dots]$ supports only
two ergodic $F$-invariant measures, namely, $\delta_Q$ and
$\delta_P$. Thus, if $\mu$ is an ergodic measure different from
$\delta_Q$ such that  $\mu([0,0,\dots])=1$, $\mu$ must be
$\delta_P$.  For this measure the Lyapunov exponent in the
central direction equals -1.

  Let us assume that $\mu([0,0,\dots])<1$. We claim that, if we define the $k$-block $\theta_k=[0,0,\dots,0,1]$, then  there exist $\eps>0$ and $l\in\N$ such that  $\mu([\theta_l])>\eps$. Indeed, just observe that
  $$
  \Lambda\setminus\{[0,0,\dots]\} = \bigcup\limits_{k=1}^\infty [\theta_k].
  $$
Thus,  there exist some positive $\eps$ a $l$-block $\theta=[0,0,\dots,0,1]$ such
that $\mu([\theta])>\eps$. By ergodicity, there exists a set
of full $\mu$-measure, $B_1\subset \Lambda$, such that every $p\in B_1$ has
 frequency for $\theta$  equal to $\mu([\theta])>\eps>0$. On the other hand,
since $\mu$ is ergodic, there exists a set $B_2\subset \Lambda$ with full $\mu$-measure such
that  for every $p\in B_2$,  the central  Lyapunov exponent is
well-defined and coincides with  $\lambda^c_\mu$.
Taking any $p \in B_1\cap B_2$ and observing
Proposition~\ref{prop.freq}, we have that
$\lambda_\mu^c<c(\mu([\theta]),l)<0$.
\end{proof}

Let us now prove item \emph{2} in Theorem \ref{t.A}.
\begin{proposition}\label{p.zero} Let $(\mu_k)$ be a sequence of ergodic measures
such that the sequence of central Lyapunov exponents $(\lambda^c_{\mu_k})$ converges to zero. Then, the sequence of
measures $(\mu_k)$  converges to $\Delta=(1/2)
\delta_Q + (1/2)\delta_P$.
\end{proposition}

\begin{proof}

Given $\eps>0$ and $\theta_l=[0,\dots,0,1]$ an $l$-block, we define $E_{\eps,l}$ by:

$$
E_{\eps,l}= \{\mu\mbox{ ergodic and $F$-invariant}; \mu(\theta_l)>\eps\}.
$$
From
Corollary~\ref{cor-lyapu-meas}   and Proposition \ref{prop.freq}, there exists a constant $a=a(l)\in(0,1)$, such that 
$$\lambda^{c}_{\mu_{ k}}=\mu_{k}([\theta_{l}])\log a.$$
Therefore $\lim_{k\rightarrow+\8}\mu_{k}([\theta_{l}])=0$. 
Since $[\theta_l]$ is open and closed  in $\Lambda$, if $\mu$ is any accumulation point for the weak* topology, we get $\mu([\theta_{l}])=0$; this holds for every $l$, which means that $\mu([\theta_l])=0$ for any $l\in \N$. Hence,  $\mu([0,0,\dots])=1$, thus $\mu=\alpha \delta_Q+(1-\alpha)\delta_P$, for
some $\alpha\in I$.

Finally, we observe that $\log|DF|_{E^{c}}|$ is continuous;  since $\mu$ is a weak* accumulation point for the sequence  $(\mu_{k})$, and $\lim_{k\rightarrow+\8}\lambda^{c}_{\mu_{k}}=0$, we get 
$$
0 =  \int \log|DF|_{E^{c}}| \,d\mu = \alpha\lambda_{\delta_Q}^c
 + (1-\alpha)\lambda_{\delta_P}^c.
$$
Since $\lambda_{\delta_Q}^c=1$ and $\lambda_{\delta_P}^c=-1$,
we must have $\alpha=1/2$. In particular, the sequence $(\mu_{k})$ admits a unique accumulation point for the weak* topology. It thus converges to $\Delta$ and  the
proof is finished.
\end{proof}

\begin{remark}
Using the structure provided by the heteroclinic cycle and the explicit  expression of $F$, we
 can prove that there exists a sequence of periodic points $p_n$ such that the Lyapunov exponents of the sequence of measures
 $\mu_n=(1/n) \sum_{i=0}^{n-1} \delta_{f^i(p_n)}$ converges to zero.
\end{remark}


\section{Proofs of theorems \ref{t.A'} and \ref{t.C}}

\subsection{Existence of equilibrium states}
In this section we prove that the entropy function $\mu \rightarrow
h_{\mu}(F)$ is upper-semicontinuous. As a consequence, we are able
to prove the existence of equilibrium states for any continuous
potential.

Observe that $F$ is not a expansive map. It can be easily deduced
observing that  points in the central segment connecting $Q$ and $P$
have same $\alpha$ and $\omega$ limits, and $F$ (respectively,
$F^{-1}$) is a contraction when it is restricted to a neighborhood of $P$
(respectively, $Q$).  Nevertheless, we have the following:

\begin{lemma}\label{lem-parti-gerador}
 Let $\mu$ be any $F$-invariant probability. Then every partition $\CP$ with diameter smaller than $1/2$ is
generating for $\mu$.
\end{lemma}
\begin{proof}
Let $\CP$ be any partition with diameter smaller
than $1/2$. For any $x$ in $\Lambda$ we denote by $\CP(x)$ the
unique element of the partition which contains $x$.
If $n$ and $m$ are two positive integers, we set
$$\CP_{-m+1}^{n-1}(x):=\bigcap_{k=-m+1}^{n-1}F^{-k}(\CP(F^k(x))),$$
and $\disp \CP_{-\8}^{+\8}(x)$ is the intersection of all $\disp \CP_{-m}^{n}(x)$.
We have just to prove that for $\mu$ almost every point $x$, $\disp\CP_{-\8}^{+\8}(x)=\{x\}$.

Consider the set of recurrent points in $\Lambda$ for $F$.
This set has full $\mu$-measure. Moreover, if $x$ is recurrent, then
its projection $\Pi(x)$ in $\Sigma_{11}$ is also recurrent. If the
bi-infinite sequence $\rho(x)$ contains at least one  $1$, Proposition
\ref{prop-renaud} proves that $\disp \CP_{-\8}^{+\8}(x)\cap
W^c(x)=\{x\}$. Hence, the uniform hyperbolicity in the two other
directions yields $\disp \CP_{-\8}^{+\8}(x)=\{x\}$.

If the bi-infinite sequence $\Pi(x)$ does not contain any $1$, then
$x$ must be in the segment $[Q,P]$. Therefore, $x=P$ or $x=Q$. Let
us first assume that $x=Q$; then for any $y\in (Q,P]$, $\lim_\ninf
F^n(y)=P$. Hence,
$$
\disp\bigcap_{n\geq
0}F^{-n}(\CP(Q))\cap[Q,P]=\{Q\}.
$$

 Again, the uniform hyperbolicity
in the two other directions yields $\disp \CP_{-\8}^{+\8}(Q)=\{Q\}$.
If $x=P$, then for any $y\in [Q,P)$, $\lim_\ninf F^{-n}(y)=Q$. The
same argument yields $\disp \CP_{-\8}^{+\8}(P)=\{P\}$.
\end{proof}

Following Proposition 2.19 in \cite{BO75}, we deduce that the metric
entropy is a upper-semicontinuous function defined on a compact set. Thus, it
attains its maximum. This imply the
existence of equilibrium states for any continuous potential and
uniqueness for any potential in a residual set of $C^0(M)$ is a standard matter, since $(\phi,\mu)\rightarrow h_\mu(f)+\int \phi \,d\mu$ is upper-semicontinuous on the set of invariant measures  and is a convex function for $\phi \in C^0(M)$.

\subsection{Phase transition: proof of Theorem~\ref{t.C}}

We denote by $\CP(t)$ the topological pressure of $\phi_{t}=t\log|DF|_{E^c}|$. For convenience it is also referred as the topological $t$-pressure. 

The function $t\mapsto \CP(t)$ is convex, thus continuous on $\R$. Hence we can define $t_0\le +\8$ as the supremum of the set 
$$\CT=\{\xi> 0,\ \forall\ t\in [0,\xi),\ \CP(t)>t\}.$$
By continuity the set $\CT$ is not empty because $\CP(0)=h_{top}(F)>0$.

\begin{lemma}\label{lem1-deltaqpas-equ}
For  $t$ in $[0,t_{0})$, any  equilibrium state $\mu_t$ for $\phi_t$ is singular with respect to $\delta_{Q}$. 
\end{lemma}
\begin{proof}
Let us assume, by contradiction, that $\mu_{t}$ is an equilibrium state for $\phi_t$ with $\mu_{t}(\{Q\})>0$, for some $t\in[0,t_{0})$. By the theorem of decomposition of measures, there exists a $F$-invariant measure $\nu$, singular with respect to $\delta_Q$ such that $\mu_t=\mu_{t}(\{Q\})\delta_{Q}+(1-\mu_{t}(\{Q\}))\nu$. Since the metric entropy is affine, we have
\begin{eqnarray*}
t<\CP(t)&=&\mu_{t}(\{Q\})t+ \big(1-\mu_{t}(\{Q\})\big) \left(h_{\nu}(F)+\int\phi_{t}\,d\nu\right)\\
&<&\mu_{t}(\{Q\})\CP(t)+ (1-\mu_{t}(\{Q\})) \left(h_{\nu}(F)+\int\phi_{t}\,d\nu\right).
\end{eqnarray*}
In particular we get $\CP(t)<h_{\nu}(F)+\int\phi_{t}\,d\nu$, which is absurd. 
\end{proof}

\begin{cor}\label{cor1-deltaqpas-equ}
Given $t$ in $[0,t_{0})$ and $\mu_t$ any equilibrium state for  $\phi_t$,  
$$\lambda^c_{\mu_{t}}=\int\log|DF|_{E^c}|d\mu_{t}< 0.$$
\end{cor}
\begin{proof}
Let $\big (\nu_{t,\xi}\big)_{ \xi \in \Lambda}$, be the ergodic decomposition of $\mu_{t}$. Since $\mu_t(\{Q\})=0$, we have that for $\mu_t$-almost every $\xi\in \Lambda$, $\nu_{t,\xi}(\{Q\})=0$.  Corollary \ref{cor-lyapu-meas} says that for each of such $\xi$, we have $\disp \int\log|DF|_{E^c}|d\nu_{t,\xi}< 0$. 
Therefore 
$$\int\log|DF|_{E^c}|d\mu_{t}=\int_{\Lambda}(\int\log|DF|_{E^c}|\,d\nu_{t,\xi})\,d\mu_t(\xi)<0.$$
\end{proof}

\begin{lemma}\label{lem1-press-decroiss}
The function $\CP$ is decreasing on $[0,t_{0})$. 
\end{lemma}
\begin{proof}
Let $t<t'$ be in $[0,t_{0})$. Let us consider two equilibrium states for $\phi_t$ and $\phi_{t'}$, $\mu_{t}$ and $\mu_{t'}$. Then we have 
\begin{eqnarray*}
\CP(t')&=&h_{\mu_{t'}}(F)+t'\lambda^c_{\mu_{t'}}\\
&=&h_{\mu_{t'}}(F)+t\lambda^c_{\mu_{t'}}+(t'-t)\lambda^c_{\mu_{t'}}\\
&\leq&\CP(t)+(t'-t)\lambda^c_{\mu_{t'}}\\
&<&\CP(t),
\end{eqnarray*}
where the last inequality yields from Corollary \ref{cor1-deltaqpas-equ}.
\end{proof}

\begin{figure}[htbp]
\begin{center}
\includegraphics[height=2in]{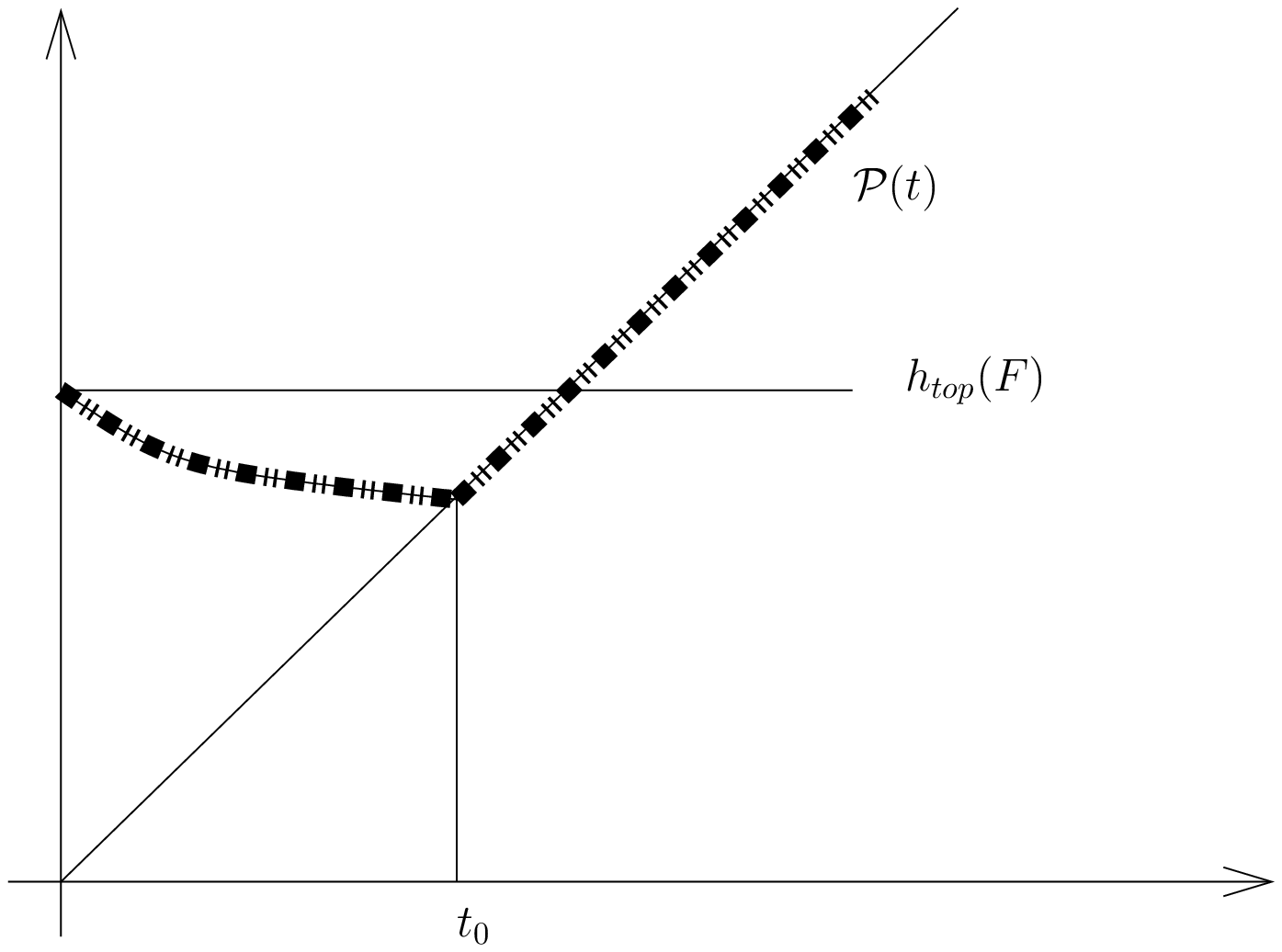}
\caption{$t\mapsto\CP(t)$}
\label{cocorico}
\end{center}
\end{figure}

Lemma \ref{lem1-press-decroiss} implies that $\CP(t)$ is less than $h_{top}(F)$ on $[0,t_{0})$. On the other hand, observe that $h_{\delta_Q}(F) + \int \phi_t \, d\delta_Q = t$, which means that $\CP(t)$ is greater or equal to $t$. Therefore $t_{0}\leq  h_{top}(F)<+\8$ (see figure \ref{cocorico} for $t\le t_0$).

We can now finish the proof of Theorem \ref{t.C}. 
Note that existence of the real number $t_{0}$ and item \emph{2} are already proved. By definition of $t_{0}$ and by continuity of $t\mapsto \CP(t)$, we must have $\CP(t_{0})=t_{0}$, thus $\delta_{Q}$ is an equilibrium state for $t_{0}$. 
Moreover, any weak$^{\*}$ accumulation point for $\mu_t$, as $t$ \emph{increases} to $t_0$, is an equilibrium state for $t_0$. Again, the continuity of $\log\left|DF|_{E^c}\right|$ yields that for such an accumulation point $\mu$, the Lyapunov exponent $\lambda^c_\mu$ is non-positive, thus the measure is different from $\delta_Q$. 

Let us pick $t>t_{0}$. 
Let $\mu_{t}$ be any equilibrium state for $t$. We have 
\begin{eqnarray*}
t\le\CP(t)&=&h_{\mu_{t}}(F)+t\lambda^c_{\mu_{t}}\\
&\le&h_{\mu_{t}}(F)+t\lambda^c_{\mu_{t}}+(t-t_{0})\lambda^c_{\mu_{t}}\\
t_0+(t-t_0)&\leq&t_{0}+(t-t_{0})\lambda^c_{\mu_{t}}.
\end{eqnarray*}
This yields  $\lambda^c_{\mu_{t}}\ge1$. Again, considering the ergodic decomposition of $\mu_{t}$, $(\nu_{t,\xi})$, we prove like in the proof of Corollary \ref{cor1-deltaqpas-equ} that for almost every $\xi$, $\nu_{t,\xi}=\delta_{Q}$. In particular, this means that $\delta_{Q}$ is the unique equilibrium state for $t>t_{0}$ (see figure \ref{cocorico} for $t\ge t_0$). This complete the proof of Theorem \ref{t.C}.


\bibliographystyle{plain}
\bibliography{bib}

\end{document}